\documentclass[12pt]{article}

\catcode`\@=11

\thicklines
\newskip\Einheit \Einheit=.6cm
\newcount\xcoord \newcount\ycoord
\newdimen\xdim \newdimen\ydim \newdimen\PfadD@cke \newdimen\Pfadd@cke
\def\PfadDicke#1{\PfadD@cke#1 \divide\PfadD@cke by2 
\Pfadd@cke\PfadD@cke \multiply\PfadD@cke by2}
\long\def\LOOP#1\REPEAT{\def\BODY{#1}\ITERATE}
\def\ITERATE{\BODY \let\next\ITERATE \else\let\next\relax\fi \next}
\let\REPEAT=\fi
\def\Punkt{\hbox{\raise-2pt\hbox to0pt{\hss\scriptsize$\bullet$\hss}}}

\def\DuennPunkt(#1,#2){\unskip
  \raise#2 \Einheit\hbox to0pt{\hskip#1 \Einheit
          \raise-1.5pt\hbox to0pt{\hss\tiny$\bullet$\hss}\hss}}
		  
\def\NormalPunkt(#1,#2){\unskip
  \raise#2 \Einheit\hbox to0pt{\hskip#1 \Einheit
          \raise-3pt\hbox to0pt{\hss\large$\bullet$\hss}\hss}}
\def\DickPunkt(#1,#2){\unskip
  \raise#2 \Einheit\hbox to0pt{\hskip#1 \Einheit
          \raise-4pt\hbox to0pt{\hss\Large$\bullet$\hss}\hss}}
\def\Kreis(#1,#2){\unskip
  \raise#2 \Einheit\hbox to0pt{\hskip#1 \Einheit
          \raise-4pt\hbox to0pt{\hss\Large$\circ$\hss}\hss}}
\def\Diagonale(#1,#2)#3{\unskip\leavevmode
  \xcoord#1\relax \ycoord#2\relax
      \raise\ycoord \Einheit\hbox to0pt{\hskip\xcoord \Einheit
         \unitlength\Einheit
         \line(1,1){#3}\hss}}
\def\AntiDiagonale(#1,#2)#3{\unskip\leavevmode
  \xcoord#1\relax \ycoord#2\relax \advance\xcoord by -0.05\relax
      \raise\ycoord \Einheit\hbox to0pt{\hskip\xcoord \Einheit
         \unitlength\Einheit
         \line(1,-1){#3}\hss}}
\def\Pfad(#1,#2),#3\endPfad{\unskip\leavevmode
  \xcoord#1 \ycoord#2 \thicklines\ZeichnePfad#3\endPfad\thinlines}
\def\ZeichnePfad#1{\ifx#1\endPfad\let\next\relax
  \else\let\next\ZeichnePfad
    \ifnum#1=1
      \raise\ycoord \Einheit\hbox to0pt{\hskip\xcoord \Einheit
         \vrule height\Pfadd@cke width1 \Einheit depth\Pfadd@cke\hss}%
      \advance\xcoord by 1
    \else\ifnum#1=2
      \raise\ycoord \Einheit\hbox to0pt{\hskip\xcoord \Einheit
        \hbox{\hskip-\PfadD@cke\vrule height1 \Einheit 
width\PfadD@cke depth0pt}\hss}%
      \advance\ycoord by 1
    \else\ifnum#1=3
      \raise\ycoord \Einheit\hbox to0pt{\hskip\xcoord \Einheit
         \unitlength\Einheit
         \line(1,1){1}\hss}
      \advance\xcoord by 1
      \advance\ycoord by 1
    \else\ifnum#1=4
      \raise\ycoord \Einheit\hbox to0pt{\hskip\xcoord \Einheit
         \unitlength\Einheit
         \line(1,-1){1}\hss}
      \advance\xcoord by 1
      \advance\ycoord by -1
    \fi\fi\fi\fi
  \fi\next}
\def\hSSchritt{\leavevmode\raise-.4pt\hbox 
to0pt{\hss.\hss}\hskip.2\Einheit
  \raise-.4pt\hbox to0pt{\hss.\hss}\hskip.2\Einheit
  \raise-.4pt\hbox to0pt{\hss.\hss}\hskip.2\Einheit
  \raise-.4pt\hbox to0pt{\hss.\hss}\hskip.2\Einheit
  \raise-.4pt\hbox to0pt{\hss.\hss}\hskip.2\Einheit}
\def\vSSchritt{\vbox{\baselineskip.2\Einheit\lineskiplimit0pt
\hbox{.}\hbox{.}\hbox{.}\hbox{.}\hbox{.}}}
\def\DSSchritt{\leavevmode\raise-.4pt\hbox to0pt{%
  \hbox to0pt{\hss.\hss}\hskip.2\Einheit
  \raise.2\Einheit\hbox to0pt{\hss.\hss}\hskip.2\Einheit
  \raise.4\Einheit\hbox to0pt{\hss.\hss}\hskip.2\Einheit
  \raise.6\Einheit\hbox to0pt{\hss.\hss}\hskip.2\Einheit
  \raise.8\Einheit\hbox to0pt{\hss.\hss}\hss}}
\def\dSSchritt{\leavevmode\raise-.4pt\hbox to0pt{%
  \hbox to0pt{\hss.\hss}\hskip.2\Einheit
  \raise-.2\Einheit\hbox to0pt{\hss.\hss}\hskip.2\Einheit
  \raise-.4\Einheit\hbox to0pt{\hss.\hss}\hskip.2\Einheit
  \raise-.6\Einheit\hbox to0pt{\hss.\hss}\hskip.2\Einheit
  \raise-.8\Einheit\hbox to0pt{\hss.\hss}\hss}}
\def\SPfad(#1,#2),#3\endSPfad{\unskip\leavevmode
  \xcoord#1 \ycoord#2 \ZeichneSPfad#3\endSPfad}
\def\ZeichneSPfad#1{\ifx#1\endSPfad\let\next\relax
  \else\let\next\ZeichneSPfad
    \ifnum#1=1
      \raise\ycoord \Einheit\hbox to0pt{\hskip\xcoord \Einheit
         \hSSchritt\hss}%
      \advance\xcoord by 1
    \else\ifnum#1=2
      \raise\ycoord \Einheit\hbox to0pt{\hskip\xcoord \Einheit
        \hbox{\hskip-2pt \vSSchritt}\hss}%
      \advance\ycoord by 1
    \else\ifnum#1=3
      \raise\ycoord \Einheit\hbox to0pt{\hskip\xcoord \Einheit
         \DSSchritt\hss}
      \advance\xcoord by 1
      \advance\ycoord by 1
    \else\ifnum#1=4
      \raise\ycoord \Einheit\hbox to0pt{\hskip\xcoord \Einheit
         \dSSchritt\hss}
      \advance\xcoord by 1
      \advance\ycoord by -1
    \fi\fi\fi\fi
  \fi\next}
\def\Koordinatenachsen(#1,#2){\unskip
 \hbox to0pt{\hskip-.5pt\vrule height#2 \Einheit width.5pt depth1 
\Einheit}%
 \hbox to0pt{\hskip-1 \Einheit \xcoord#1 \advance\xcoord by1
    \vrule height0.25pt width\xcoord \Einheit depth0.25pt\hss}}
\def\Koordinatenachsen(#1,#2)(#3,#4){\unskip
 \hbox to0pt{\hskip-.5pt \ycoord-#4 \advance\ycoord by1
    \vrule height#2 \Einheit width.5pt depth\ycoord \Einheit}%
 \hbox to0pt{\hskip-1 \Einheit \hskip#3\Einheit 
    \xcoord#1 \advance\xcoord by1 \advance\xcoord by-#3 
    \vrule height0.25pt width\xcoord \Einheit depth0.25pt\hss}}
\def\Gitter(#1,#2){\unskip \xcoord0 \ycoord0 \leavevmode
  \LOOP\ifnum\ycoord<#2
    \loop\ifnum\xcoord<#1
      \raise\ycoord \Einheit\hbox to0pt{\hskip\xcoord 
\Einheit\Punkt\hss}%
      \advance\xcoord by1
    \repeat
    \xcoord0
    \advance\ycoord by1
  \REPEAT}
\def\Gitter(#1,#2)(#3,#4){\unskip \xcoord#3 \ycoord#4 \leavevmode
  \LOOP\ifnum\ycoord<#2
    \loop\ifnum\xcoord<#1
      \raise\ycoord \Einheit\hbox to0pt{\hskip\xcoord 
\Einheit\Punkt\hss}%
      \advance\xcoord by1
    \repeat
    \xcoord#3
    \advance\ycoord by1
  \REPEAT}
\def\Label#1#2(#3,#4){\unskip \xdim#3 \Einheit \ydim#4 \Einheit
  \def\lo{\advance\xdim by-.5 \Einheit \advance\ydim by.5 \Einheit}%
  \def\llo{\advance\xdim by-.25cm \advance\ydim by.5 \Einheit}%
  \def\loo{\advance\xdim by-.5 \Einheit \advance\ydim by.25cm}%
  \def\o{\advance\ydim by.25cm}%
  \def\ro{\advance\xdim by.5 \Einheit \advance\ydim by.5 \Einheit}%
  \def\rro{\advance\xdim by.25cm \advance\ydim by.5 \Einheit}%
  \def\roo{\advance\xdim by.5 \Einheit \advance\ydim by.25cm}%
  \def\l{\advance\xdim by-.30cm}%
  \def\r{\advance\xdim by.30cm}%
  \def\lu{\advance\xdim by-.5 \Einheit \advance\ydim by-.6 \Einheit}%
  \def\llu{\advance\xdim by-.25cm \advance\ydim by-.6 \Einheit}%
  \def\luu{\advance\xdim by-.5 \Einheit \advance\ydim by-.30cm}%
  \def\u{\advance\ydim by-.30cm}%
  \def\ru{\advance\xdim by.5 \Einheit \advance\ydim by-.6 \Einheit}%
  \def\rru{\advance\xdim by.25cm \advance\ydim by-.6 \Einheit}%
  \def\ruu{\advance\xdim by.5 \Einheit \advance\ydim by-.30cm}%
  #1\raise\ydim\hbox to0pt{\hskip\xdim
     \vbox to0pt{\vss\hbox to0pt{\hss$#2$\hss}\vss}\hss}%
}
\catcode`\@=12

\usepackage[colorlinks=true,
linkcolor=green,
filecolor=brown,
citecolor=green]{hyperref}
\usepackage{amsmath}
\def\t{Tak\'{a}cs}

\newcommand{\seqnum}[1]{\href{http://www.research.att.com/cgi-bin/access.cgi/as/~njas/sequences/eisA.cgi?Anum=#1}{\underline{#1}}}

\begin{document}

\begin{center}
{\Large
 A Combinatorial Derivation of the Number of Labeled Forests      \\ 
}
\vspace{10mm}
DAVID CALLAN  \\
Department of Statistics  \\
University of Wisconsin-Madison  \\
1210 W. Dayton St   \\
Madison, WI \ 53706-1693  \\
{\bf callan@stat.wisc.edu}  \\
\vspace{5mm}

\end{center}
\begin{abstract}
    Lajos \t\  gave a somewhat formidable alternating sum expression    
	for the number of forests of unrooted trees on $n$ labeled vertices. 
	Here we use a weight-reversing involution on suitable tree configurations 
	to give a combinatorial derivation of \t' result. 
\end{abstract}

\vspace{5mm}

\t\ \cite{takacs} used Lagrange inversion to obtain the alternating sum expression
\begin{equation}
\frac{n!}{n+1} \sum_{j=0}^{\lfloor n/2 \rfloor} (-1)^{j}
\frac{(2j+1)(n+1)^{n-2j}}{2^{j} j! (n-2j)!}
\label{eq:1}
\end{equation}
for the number of forests of unrooted trees on $[n]=\{1,2,\ldots,n\}$
\htmladdnormallink{A001858}{http://www.research.att.com:80/cgi-bin/access.cgi/as/njas/sequences/eisA.cgi?Anum=A001858}.
This contrasts with Cayley's simple 
expression $(n+1)^{n-1}$
\htmladdnormallink{A000272}{http://www.research.att.com:80/cgi-bin/access.cgi/as/njas/sequences/eisA.cgi?Anum=A000272}
for the number of forests of rooted trees on $[n]$.  
Here we use well-known 
counts for forests of rooted trees to give a combinatorial derivation 
of \t's result:
we present (\ref{eq:1}) as the total weight of certain 
weighted tree configurations in which forests of unrooted trees show up with 
weight $+1$ and we exhibit a weight-reversing
involution that cancels out the weights of all other configurations.
First, rewrite (\ref{eq:1}) as 
\begin{equation}
 \sum_{0\le j \le 
n/2}^{}\ 
(-1)^{j}\phantom{\int_{}^{}}\underbrace{\binom{n}{2j}(2j-1)!!}_{A}\phantom{\int_{}^{}}
\underbrace{(2j+1)(n+1)^{(n+1)-(2j+1)-1}\phantom{\int_{}^{}}}_{B}
	\label{eq:2}
\end{equation}
where $(2j-1)!!=1\cdot 3\cdot 5  \ldots  (2j-1)$ is the 
double factorial. The factor $B$ is the number of forests on $[0,n]$ 
consisting of $2j+1$ trees rooted at a specified set of $2j+1$ roots 
\cite[Theorem 3.3, p.\:17]{moon}
(see also \cite[\S 2.1]{pitman} for a recent elegant proof).
The factor $A$ is the number of ways to select $2j$ elements 
from $[n]$ and then divide them up into pairs; in other words, 
to form a perfect matching on some $2j$ elements of $[n]$.
These $2j$ elements, together with 0, serve nicely as the specified 
roots to construct configurations counted by the product $AB$.

Define a \emph{partially-paired rooted} (PPR, for short) \emph{$n$-forest}  
to be a tree rooted at 0 
and zero or more (unordered) pairs of rooted trees, the vertex sets of all the trees 
forming a partition of $[0,n]$. 
The \emph{pair-count} of a PPR forest is its number of pairs of trees.
The product $AB$ is then the number of PPR $n$-forests with pair-count
$j$.  If we define the \emph{weight} of a PPR forest
of pair-count $j$ to be $(-1)^{j}$, then the right hand side of (1) is the
total weight of all PPR $n$-forests.

To include the objects we're trying to count
among these PPR $n$-forests, we suppose each tree in an 
unrooted forest to be rooted at its \emph{smallest} vertex.
Then forests of unrooted trees on $[n]$ correspond precisely to
PPR $n$-forests with pair-count 0 and each 
child of vertex 0 smaller than all its descendants (delete vertex 0 to get the 
forest of unrooted trees). A vertex $v$ in a 
rooted tree is \emph{inversion-initiating} if at least one descendant of $v$ is 
$<v$, otherwise it is \emph{regular}.
Thus forests of unrooted trees on $[n]$ appear as
PPR $n$-forests with pair-count 0 and all children of vertex 0 regular.
These special PPR forests are counted with weight 1 
and here is the promised weight-reversing involution on the rest.

Given a PPR forest, let $a$ denote the smallest vertex among all trees (if any) other than the one 
rooted at 0, let $u$ be the root of $a$'s tree ($u$ is possibly, 
but not necessarily, = $a$), and let $v$ be the root of the other tree in its pair.
At the same time, if 0 has any  inversion-initiating
children, let $a'$ be the smallest among 
all descendants of these inversion-initiating
vertices, let $v'$ be the child of 0 of which $a'$ is a descendant, 
and let $u'$ (possibly = $a'$) be the child of $v'$ on the path from 
$v'$ to $a'$. See the illustration below where solid lines represent 
mandatory edges, vertical dotted lines optional edges, and diagonal 
dotted lines optional subtrees. 


\vskip30pt
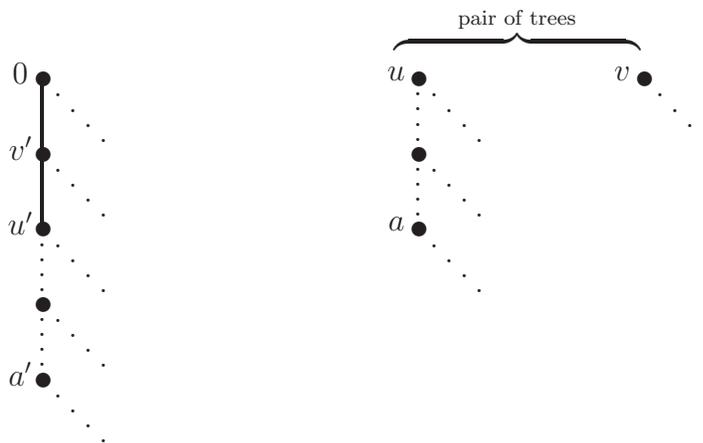
\vbox{
$$
\Einheit=1cm.
\PfadDicke{1pt}
\SPfad(-5,-4),22\endSPfad
\SPfad(-5,-1),4\endSPfad
\SPfad(-5,-2),4\endSPfad
\SPfad(-5,-3),4\endSPfad
\SPfad(-5,-4),4\endSPfad
\SPfad(-5,0),4\endSPfad
\SPfad(0,-2),22\endSPfad
\SPfad(0,0),4\endSPfad
\SPfad(0,-1),4\endSPfad
\SPfad(0,-2),4\endSPfad
\SPfad(3,0),4\endSPfad
\Pfad(-5,-2),22\endPfad
\Label\r{\overbrace{\phantom{\text{pair of trees yeah }}}^{\text{pair of trees}}}(1,.5)
\Label\l{0}(-5,0.1) 
\Label\l{u}(0,0.1)
\Label\l{v}(3,0.1)
\Label\l{v'}(-5,-.9) 
\Label\l{u'}(-5,-1.9) 
\Label\l{a'}(-5,-3.9) 
\Label\l{a}(0,-1.9) 
\NormalPunkt(-5,0)
\NormalPunkt(-5,-1)
\NormalPunkt(-5,-2)
\NormalPunkt(-5,-3)
\NormalPunkt(-5,-4)
\NormalPunkt(0,0)
\NormalPunkt(0,-1)
\NormalPunkt(0,-2)
\NormalPunkt(3,0)
\Label\o{\textrm{$a'$ is smallest descendant of an}}(-5,-5.7) 
\Label\o{\textrm{inversion-initiating child of 0 }}(-5,-6.2) 
\Label\o{\textrm{$a$ is smallest vertex not}}(1,-5.7) 
\Label\o{\textrm{a descendant of 0}}(1,-6.2) 
$$
}
\vskip5pt

At least one of $a,a'$ will exist unless the pair-count is 0 and all 
children of vertex 0 are regular; these are the special PPR forests, 
representing unrooted forests, and they survive.
Choose the smaller of $a,a'$. If it's $a$, add an edge from $0$ to $v$ 
and an edge
from $v$ to $u$ so that vertex 0 acquires a new inversion-initiating child $v$ (with a 
small descendant $a$) and the number of pairs of trees is 
reduced by 1. If it's $a'$, delete the edges $0v'$ and 
$v'u'$ to form a new pair of trees rooted at $u'$ and $v'$ (with $a'$ 
now the smallest vertex among all pairs of trees). In either case, 
the number of pairs of trees changes by 1, so the weight changes 
sign. The mapping is clearly an involution on all non-special 
PPR forests and so their weights cancel out. 
Thus (\ref{eq:2}) (= (\ref{eq:1})) is the  number of forests of unrooted trees on $[n]$.

\bigskip
\hrule
\bigskip



\bigskip
\hrule
\bigskip

\noindent (Concerned with sequence \seqnum{A001858}.)

\end{document}